\documentclass[11pt, a4paper]{amsart}
\usepackage{amsfonts,amssymb,amsmath,amsthm,amscd,mathtools,multicol,tikz, tikz-cd,caption,enumerate,mathrsfs,thmtools}
\usepackage{inputenc}
\usepackage[foot]{amsaddr}
\usepackage[doi=false,isbn=false,url=false]{biblatex} 
\usepackage[ colorlinks, linkcolor=blue, citecolor=red]{hyperref}
\usepackage{latexsym}
\usepackage{fullpage}
\usepackage{microtype}
\usepackage{subfiles} 
\usepackage[lined,boxed,commentsnumbered]{algorithm2e}
\usepackage{booktabs}
\usepackage{xcolor}

\counterwithin{table}{section}
\DeclareMathOperator{\SL}{SL}
\DeclareMathOperator{\GL}{GL}
\DeclareMathOperator{\PSL}{PSL}
\DeclareMathOperator{\PGL}{PGL}
\DeclareMathOperator{\PSU}{PSU}

\usepackage[T1]{fontenc}

\usepackage{palatino}
\makeatletter
\makeindex 

\newcommand{\be}{\begin{equation}}
\newcommand{\ee}{\end{equation}}
\newcommand{\beano}{\begin{eqn*}} 
	\newcommand{\eeano}{\end{eqnarray*}}
\newcommand{\ba}{\begin{array}}
	\newcommand{\ea}{\end{array}}

\declaretheoremstyle[headfont=\normalfont]{normalhead}

\newtheorem{theorem}{Theorem}[section]

\newtheorem{lemma}[theorem]{Lemma}

\newtheorem{definition}[theorem]{Definition}

\newcommand{\Inn}{\mathrm{Inn}}

\numberwithin{equation}{section}

\newtheorem{ex}[theorem]{Example}

\newtheorem{problem}[theorem]{Problem}
\@namedef{subjclassname@2020}{%
  \textup{2020} Mathematics Subject Classification}
\makeatother

\begin{document}
\title{Classification of simple quandles of small order}
\author{Dilpreet Kaur}

\author{Pushpendra Singh}

\subjclass[2020]{20N02,20B35,57K12,20D05}
\keywords{Quandles, simple quandles, simple groups, quasiprimitive groups}
\newcommand{\Addresses}{{
  \bigskip
  \footnotesize

  ( Dilpreet Kaur ) \textsc{Department of Mathematics, Indian Institute of Technology Jodhpur, Jodhpur, Rajasthan, 342030, India}\par\nopagebreak
  \textit{E-mail address}: \texttt{dilpreetkaur@iitj.ac.in}

  \medskip

  ( Pushpendra Singh ) \textsc{Department of Mathematics, Indian Institute of Technology Jodhpur, Jodhpur, Rajasthan, 342030, India}\par\nopagebreak
  \textit{E-mail address}: \texttt{singh.105@iitj.ac.in}

}}
\maketitle

\begin{abstract}
In this article, we define quasiprimitive quandles and provide a description of simple quandles with the help of quasiprimitive groups. Consequently, we enumerate finite non affine imprimitive simple quandles up to order $4096$. We then describe properties of simple quandles arising from $\PSL(2,q)$ and $\PGL(2,q)$.
\end{abstract}

\section{Introduction}\label{1}
Quandles are non associative algebraic structures derived from the axiomatization of Reidemeister moves of planar diagrams of knots and links. They are introduced independently in \cite{Joy82a} and \cite{Mat82}. Apart from their applications in knot theory for constructing invariants of knots and links, they are rich algebraic objects in their own right. Simple quandles are introduced in \cite{Joy82b} with the definition motivated from universal algebra. A quandle is simple if it has no nontrivial congruences. Finite simple quandles are building blocks of general quandles \cite{BS21}. Quandles have deep connections with group theory. A general description of simple quandles is given in \cite{Joy82b} and \cite{AG03} with the help of finite simple groups. 

In this article, we describe simple quandles with the help of quasiprimitive permutation groups and prove that simple quandles are quasiprimitive. Further, we enumerate finite non affine imprimitive simple quandles up to order 4096. We use a result from \cite{KS25} to construct multiplication tables efficiently and provide a novel, computationally efficient isomorphism theorem, which allows us to complete the enumeration in a few hours. 

We observe from the data that many simple quandles have $\PSL(2,q)$ and $\PGL(2,q)$ as inner automorphism groups. We identify conjugacy classes in these groups that correspond to primitive quandles. Furthermore, we describe isomorphism classes and \textit{profile} of these quandles.

\section{Review of quandles } \label{2}
A quandle is an algebraic structure $(X,\triangleright)$ satisfying the following axioms:
\begin{enumerate}[(1)]
\item $x \triangleright x = x$ for all $x \in X$.
\item The map $R_x : X \to X, y \to y \triangleright x$ is a bijection of $X$ for all $x \in X$.
 \item $(x~\triangleright~y)\triangleright~z = (x~\triangleright~z)\triangleright(y \triangleright z)$  for all $x,y,z \in X.$
\end{enumerate}

The axioms $(1),(2)$ and $(3)$ are motivated from the Reidemeister moves \text{I}$,$\text{II} and \text{III} respectively. The axiom $(2)$ induces a binary operation ${\triangleright}^{-1}$ on $X$ such that $(X,\triangleright^{-1})$ also forms a quandle structure. This operation is defined as $x~{\triangleright}^{-1}~y = z$ if $x = z \triangleright y$. 
The structure obtained by removing axiom $(1)$ is known as a \textit{rack}. 

A quandle homomorphism is a mapping $f : X \to X$ satisfying $f(x\triangleright y) = f(x) \triangleright f(y)$ and $f(x~{\triangleright^{-1}}~y) = f(x)~{\triangleright^{-1}}~f(y)$. The mappings in axiom $(2)$ are known as inner automorphisms. They generate a subgroup $\operatorname{Inn}(X)$ of $\operatorname{Aut}(X)$, which is called the inner automorphism group.
$$\operatorname{Inn}(X) = \langle R_x~:~x \in X \rangle$$

A quandle is called connected ( resp. primitive) if the natural action of $\operatorname{Inn}(X)$ on $X$ is a transitive ( resp. primitive) group action. The displacement group $\operatorname{Dis}(X)$ of a quandle $(X,\triangleright)$ is a subgroup of $\operatorname{Inn}(X)$ generated by elements $R^{-1}_xR_y$ for $x,y$ in $X$. The quotient group $\operatorname{Inn}(X)/\operatorname{Dis}(X)$ is a cyclic group. In case of connected quandles, $\operatorname{Dis}(X)=\operatorname{Inn}(X)'$ and $\operatorname{Dis}(X)$ also acts transitively on $X$ \cite{Joy82b}. 

\begin{ex}\label{eg1}
Let $G$ be a group and $\psi \in \operatorname{Aut}(G)$. Let $H$ be the subgroup of $G$ containing elements fixed by $\psi$. Then the set of cosets $G / H$ forms a quandle with the operation $xH \triangleright yH = y\psi(y^{-1}x)H$. It is denoted by $(G, H,\psi)$.
\end{ex}

\begin{ex}\label{eg2}
Let $G$ be a group and $C$ be a subset of $G$. Let $C$ be closed under conjugation. Then $C$ forms a quandle with conjugation as a quandle operation. It is denoted by $\operatorname{Conj}(G, C)$. If $C$ is taken to be the whole group, then it is denoted by $\operatorname{Conj}(G)$.
\end{ex}

\begin{definition}
\textit{Profile} of a quandle: Let $(X,\triangleright)$ be a finite quandle. Then the set of partition types of all right translations is called \textit{profile} of $X$.
\[ {\rm Prof}(X) := \{\Lambda_x \;;\; x \in X\} \] where $\Lambda_x$ is the partition type of corresponding $R_x \in {\rm Sym}(X)$.
\end{definition}
For a connected quandle, all right translations have the same partition type. In that case, the \textit{profile} is simply $\Lambda_e$ for a fixed $e \in X$. 
If the natural homomorphism $X \to \operatorname{Conj}(\operatorname{Inn}(X)), x \mapsto R_{x}$ is injective, then $(X,\triangleright)$ is called a faithful quandle. A quandle $(X,\triangleright)$ is called simple if $|X|>1$ and every quandle homomorphism from $X$ is either injective or constant. Simple quandles of size $>2$ are connected. Simple quandles are faithful \cite{Joy82b}. In Example \ref{eg1}, if we consider $G$ to be an abelian group and $H$ to be the trivial subgroup, then the quandle $(G,\{e\},\psi)$ is known as an affine quandle.

\begin{definition}
The action of a group $G$ on a set $X$ is called a quasiprimitive group action if it is faithful and the restriction of the action to every nontrivial normal subgroup is a transitive group action.
\end{definition}

\begin{definition}
We call a quandle $(X,\triangleright)$ a quasiprimitive quandle if the action of the group $\operatorname{Inn}(X)$ on $X$ is a quasiprimitive group action.
\end{definition}

For a primitive group $G$ on a set $X$, every nontrivial normal subgroup of $G$ also acts transitively on $X$ \cite[Theorem 1.7]{Cam99}. Thus, we have the following relationship among these groups.
$$\{A_n,S_n\} \subseteq \text{Primitive Groups} \subseteq \text{Quasiprimitive Groups} \subseteq \text{Transitive Groups}$$

The following theorems characterize the $\operatorname{Inn}(X)$ when $X$ is a simple quandle.

\begin{theorem}\cite{Joy82b}\label{t2}
Let $G$ be a finite group then $G \cong \operatorname{Inn}(X)$ for some simple quandle $(X,\triangleright)$ iff the following are satisfied
\begin{enumerate}
\item
$Z(G)$ is trivial.
\item
$G/G'$ is cyclic.
\item
$G'$ is the smallest nontrivial normal subgroup of $G$.
\end{enumerate}

\end{theorem}

\begin{theorem}\cite{AG03}\label{t3}
Let $(X,\triangleright)$ be a non affine simple quandle. Then $G \cong \operatorname{Inn}(X)$ is of the form $L^t \rtimes \langle \phi   \rangle /Z(L^t \rtimes \langle \phi   \rangle)$ where $L$ is a non abelian simple group, $t \in \mathbb{N}$,$  \phi \in \operatorname{Aut}(L^t)$ and $\phi$ acts on $L^t$ by $\phi(l_1,l_2,\hdots,l_t) = (\theta(l_t),l_1,\hdots,l_{t-1})$ for $\theta \in \operatorname{Aut}(L)$. Furthermore $G'=L^t$.
\end{theorem}

\begin{theorem}\cite[Theorem 1.1.1]{Bor17} \label{bors}
Let $G$ be a finite group. If $G$ has an automorphism of order greater than $\frac{1}{2}|G|$, then $G$ is abelian.
\end{theorem}

\section{Classification of simple quandles} \label{5}

We begin with the following theorem, which allows us to use quasiprimitive groups to classify finite simple quandles.

\begin{theorem}
Let $X$ be a finite simple quandle, then ${\rm Inn}(X)$ is a quasiprimitive group.
\end{theorem}

\begin{proof}
In Theorem \ref{t3}, we note that $Z(L^t \rtimes \langle \phi \rangle)$ is isomorphic to a subgroup of $\langle \phi  \rangle$. Moreover $|\langle \phi \rangle |< \frac{1}{2}|L^t|$ using Theorem \ref{bors}. Thus, we get that for a simple quandle $(X,\triangleright)$, if $K$ is a normal subgroup of $G = \operatorname{Inn}(X)$, then $G' \leq K \leq G$. Since simple quandles are connected, this further implies that all nontrivial normal subgroups of $\operatorname{Inn}(X)$ act transitively on $X$.
\end{proof}

As a consequence of the above theorem, we get that simple quandles are quasiprimitive quandles. However, we note that there are quasiprimitive quandles that are not simple, for example, the \texttt{RIG} \cite{Ven14} quandle $(20,3)$. Quasiprimitive groups are classified in \cite{Pra93} and can be divided into two classes, namely primitive groups and quasi imprimitive groups.

Now we give the following key lemma, which will be useful to give a fast isomorphism theorem for finite simple quandles.

\begin{lemma}\label{isoqthm}
Let $(X,\triangleright)$ and $(Y,\triangleright)$ be quandles. Then a quandle isomorphism $f: X\to Y$ extends to a group isomorphism $ \phi : \operatorname{Inn}(X) \to \operatorname{Inn}(Y)$.
\end{lemma}

\begin{proof}
We have $\operatorname{Inn}(X)=\langle R_{x}:x\in X \rangle$ and $\operatorname{Inn}(Y)= \langle R_{y}: y \in Y  \rangle$. Define $\phi : \operatorname{Inn}(X) \to \operatorname{Inn}(Y)$ as $\phi(R_{x_1}^{\varepsilon_1}R_{x_2}^{\varepsilon_2}\hdots R_{x_m}^{\varepsilon_m}) = R_{f(x_1)}^{\varepsilon_1}R_{f(x_2)}^{\varepsilon_2}\hdots R_{f(x_m)}^{\varepsilon_m} $ where $x_i \in X, \varepsilon_i = \pm 1.$ Let $\triangleright^{1} = \triangleright$, ${\triangleright}^{-1}$ on $X$ is defined as $x~{\triangleright}^{-1}~y = z$ if $x = z \triangleright y$. It is enough to show that $\phi$ is a well-defined mapping. Suppose $R_{x_1}^{\varepsilon_1}R_{x_2}^{\varepsilon_2}\hdots R_{x_m}^{\varepsilon_m} = R_{y_1}^{\varepsilon{'}_1}R_{y_2}^{\varepsilon{'}_2}\hdots R_{y_k}^{\varepsilon{'}_k}$. Let $y=f(x) \in Y$. We have $((( y \triangleright^{\varepsilon_m}~f(x_m))\triangleright^{\varepsilon_{m-1}}~f(x_{m-1})) \hdots ) \triangleright^{\varepsilon_1}~f(x_1)=f((((x \triangleright^{\varepsilon_m}~x_m)\triangleright^{\varepsilon_{m-1}}~x_{m-1}) \hdots ) \triangleright^{\varepsilon_1}~x_1 ) = f((((x \triangleright^{\varepsilon{'}_k}~y_k)\triangleright^{\varepsilon{'}_{k-1}}~y_{k-1}) \hdots ) \triangleright^{\varepsilon{'}_1}~y_1 ) = (((y \triangleright^{\varepsilon{'}_k}~f(y_k))\triangleright^{\varepsilon{'}_{k-1}}~f(y_{k-1}))\hdots ) \triangleright^{\varepsilon{'}_1}~f(y_1).$ This implies $R_{f(x_1)}^{\varepsilon_1}R_{f(x_2)}^{\varepsilon_2}\hdots R_{f(x_m)}^{\varepsilon_m} = R_{f(y_1)}^{\varepsilon{'}_1} R_{f(y_2)}^{\varepsilon{'}_2}\hdots R_{f(y_k)}^{\varepsilon{'}_k}$. Hence $\phi$ is an isomorphism. 
\end{proof}

We note that \cite[Section 5]{HSV16} can be restricted to quasiprimitive quandles with the help of quasiprimitive groups. Let $X$ be a finite set with $|X|>2$ and $G$ be a finite quasiprimitive group on the set $X$. Let $G_e$ be the stabilizer group of $e \in X$. For a tuple $(G,\rho)$ where $\rho \in Z(G_e)$ with $\langle \rho^G \rangle = G$, set $X$ can be given a quandle structure in following way
$$\mathcal{QPQ}(G,\rho) = (X,\circ),~~~x \circ y = {\alpha_{y}}{\rho}{\alpha^{-1}_{y}}(x)$$
where $\alpha_{y} \in G$ such that $\alpha_y(e)=y$. Moreover ${\Inn}(X)=G$.

Let $Q=(G,G_e,\psi_\rho)$ be coset quandle with operation $gG_e \triangleright hG_e=h\rho h^{-1}g\rho^{-1}G_e$. We have ${\rm Inn}(Q)=\langle R_{gG_e}\;; gG_e\in G/G_e \rangle$. The mapping $f:(G,G_e,\psi_\rho) \rightarrow (X,\circ), gG_e \mapsto g(e)$ is a quandle isomorphism similar to \cite[Lemma 3.3]{KS25}. The mapping $\phi: {\rm Inn}(Q) \rightarrow G, R_{gG_e}\mapsto g\rho g^{-1} $ is a group isomorphism. The pair $(\phi,f)$ is a group action isomorphism between the action of ${\rm Inn}(Q)$ on $G/G_e$ and the action of $G$ on $X$.

\begin{theorem}
Let $G$ be a quasiprimitive group of degree $n$ on the set $X$. Let $\rho \in Z(G_e)$ with $\langle \rho^{G} \rangle = G$. A quandle is a quasiprimitive quandle of order $n$ if and only if it is isomorphic to $\mathcal{QPQ}(G,\rho)$.
\end{theorem} 

\begin{proof}
The proof is similar to \cite[Theorem 3.3]{KS25}.
\end{proof}

In the case when $\mathcal{QPQ}(G,\rho)$ is a faithful connected quandle, we have the following isomorphism.

\begin{lemma}\label{isoconj}
Let $(X,\circ)$ be a faithful connected quandle then $(X,\circ) \cong \operatorname{Conj}(G,\rho^G)$.
\end{lemma}

\begin{proof}
The proof is similar to Lemma \cite[Lemma 3.5]{KS25}.
\end{proof}

\begin{theorem} \label{isoq}
Let $\mathcal{QPQ}(G,\rho)$ and $\mathcal{QPQ}(G,\rho{'})$ be faithful quasiprimitive quandles. Then  $\mathcal{QPQ}(G,\rho) \cong \mathcal{QPQ}(G,\rho{'})$ if and only if there exists an automorphism $\phi$ of group $G$ satisfying $\phi(\rho^G) = \rho{'}^{G}$.
\end{theorem}

\begin{proof}
The reverse direction is straightforward. The forward result follows as the consequence of Lemma \ref{isoqthm} and Lemma \ref{isoconj}.
\end{proof}

The classification of non affine primitive quandles up to order 4096 is given in \cite{KS25}. Quasi imprimitive groups are classified up to order 4096 in \texttt{GAP} \cite{GAP} package \texttt{QuimpGrp} \cite{Ber22}. We implement the algorithm from \cite{KS25} in \texttt{GAP} on an Intel Core i$7$-$10700$ $2.90$ GHz processor. We use Theorem \ref{isoq} for isomorphism checking of faithful quasi imprimitive quandles. 

We construct faithful quasi imprimitive quandles up to order 4096. We generate 991 quandles, out of which 802 are up to isomorphism. It turns out they are all simple quandles by using Theorem \ref{t2}. The description of $\operatorname{Inn}(X)$ and $\operatorname{Dis}(X)$ is provided in Tables \ref{tables1} to \ref{tables8}. $[n,i]$ denotes $i^{th}$ quandle of order $n$.  The database of constructed quandles is available at \cite{DKa}. The complete search takes around 14 hours and 30 minutes of computation time. 

The \texttt{QuimpGrp} data can also be used to enumerate nonfaifthful connected quasiprimitive quandles. The following theorem is efficient for isomorphism checking for (nonfaithful) connected quandles.

\begin{theorem}
Let $Q_1=(G,G_e,\psi_{\rho_1})$ and $Q_2=(G,G_e,\psi_{\rho_2})$ be two connected quandles constructed from a transitive group $G$. Then $Q_1 \cong Q_2$ if and only if there is $\phi \in {\rm Aut}(G)$ such that $\phi(G_e)=G_e$ and $\phi(\rho_1)=\rho_2$. 
\end{theorem}

\begin{proof}
Suppose we are given $\phi \in {\rm Aut}(G)$ such that $\phi(G_e)=G_e$ and $\phi(\rho_1)=\rho_2$. Now define $f: (G,G_e,\psi_{\rho_1}) \rightarrow (G,G_e,\psi_{\rho_2}), gG_e \mapsto \phi(g)G_e$. We first show that $f$ is well defined. Let $gG_e=hG_e$, then $h^{-1}g \in G_e$, so $\phi(h)^{-1}\phi(g)=\phi(h^{-1}g) \in G_e$. Thus $\phi(g)G_e=\phi(h)G_e$. Its inverse is defined as $f^{-1}(gG_e)=\phi^{-1}(g)G_e$. Now we prove $f$ is a quandle homomorphism. We have $f(gG_e\triangleright hG_e)=f(h\rho_1h^{-1}g{\rho_1^{-1}}G_e)=\phi(h\rho_1h^{-1}g{\rho_1^{-1}})G_e=\phi(h)\phi(\rho_1)\phi(h^{-1})\phi(g)\phi(\rho_1^{-1})G_e=\phi(h)\rho_2\phi(h)^{-1}\phi(g)\rho_2^{-1}G_e=\phi(g)G_e\triangleright \phi(h)G_e=f(gG_e)\triangleright f(hG_e)$.

Conversely, suppose $f:(G,G_e,\psi_{\rho_1})\rightarrow (G,G_e,\psi_{\rho_2})$ be a quandle isomorphism. Then as a consequence of Lemma \ref{isoqthm}, $f$ induces a group isomorphism $\phi: {\rm Inn}(Q_1) \rightarrow {\rm Inn}(Q_2), R_{gG_e} \mapsto R_{f(gG_e)}$ and so for any $u \in {\rm Inn}(Q_1)$, we have $\phi(u)f(x)=f(u(x))$. We may assume $f(G_e)=G_e$, because if not, then since the action of $G$ on $G/G_e$ is transitive, there is $a\in G$ such that $a(f(G_e))=G_e$ and $a\circ f$ is a quandle isomorphism. Now $\phi(R_{G_e})=R_{f(G_e)}=R_{G_e}$. Now let $H_i={\rm Stab}_{\rm Inn(Q_i)}(G_e)$ and $h\in H_1 $, then $h(G_e)=G_e$. We have $\phi(h)(G_e)=\phi(h)(f(G_e))=f(h(G_e))=f(G_e)=G_e$. So $\phi(H_1) \subseteq H_2$. Similarly $H_2 \subseteq \phi(H_1)$ using $\phi^{-1}$. Thus $\phi(H_1)=H_2$. Now the result follows by identifying ${\rm Inn}(Q_i)$ with $G$.
\end{proof}

We end this section by addressing a problem posed in \cite{KS25}. The authors in \cite[Problem 4.1]{KS25} ask for which groups $G$, does the bijection $f: C_1 \to C_2$ satisfying $f(yxy^{-1})=f(y)f(x)f(y)^{-1}$ for all $x,y \in C_1$, extends to an automorphism of $G$, where $C_1,C_2$ are conjugacy classes of $G$, satisfying $\langle C_1 \rangle = \langle C_2 \rangle =G$. Let $G$ be a finite group with trivial center and $C$ be a conjugacy class of $G$ satisfying $\langle C \rangle =G$. Then for the quandle $X=\operatorname{Conj}(G,C)$, we have $\operatorname{Inn}(X) \cong G/Z(G)=G$. Thus, as a consequence of the above Lemma \ref{isoqthm}, the above mentioned problem holds true.

\section{Simple quandles arising from $\PSL(2,q)$ and $\PGL(2,q)$}

We note that many simple quandles in the data have ${\rm Inn}(X)$ equal to either $\PGL(2,q)$ or $\PSL(2,q)$. Therefore, we analyze such simple quandles in detail in this section. We begin with classifying conjugacy classes that give primitive quandles.

Let $\mathbb{F}_q$ be a finite field\index{Finite field}. Let $\GL(2,q)$ be the group of $2 \times 2$ matrices over $\mathbb{F}_q$ with nonzero determinant. Consider $\SL(2,q)$ is the subgroup of $\GL(2,q)$ consisting of determinant one matrices. Their projective counterparts are obtained by taking the quotient by their center. Thus we have $\PSL(2,q) = \SL(2,q)/ \{\pm I\}$ and $\PGL(2,q)=\GL(2,q)/\{{\lambda}I~:\lambda \in \mathbb{F}^{*}_q\}$. In this section, we take $p$ to be an odd prime and $q=p^n$ for some $n \in \mathbb{N}$. We choose a $\Delta \in \mathbb{F}_q/\mathbb{F}^2_{q}$.

We begin with classifying primitive quandles. The following lemma is helpful for this classification.

\begin{lemma}\cite[Lemma 4.1]{CS25} \label{c5s5.6}
Let $G$ be a group and $x \in G$ such that $\langle x^G \rangle=G$. If $\{x\} \subsetneq \langle x \rangle \cap x^G \subsetneq x^G$ then $\operatorname{Conj}(G,x^G)$ is imprimitive.
\end{lemma}

The conjugacy classes of $\PSL(2,q)$ are computed in [\cite{Ada02},\cite{FGV10}] and described in Table \ref{tpsl1} and \ref{tpsl2}.

\begin{table}[ht]
$$
\begin{array}{|c|c|c|c|c|} 
\hline
 & \text{Representative} & \text{No.} & \text{Size} & \text{Order representative} (d)\\
\hline
\mathcal{C}_1 &c_1 = \left ( \begin{array}{rr}
                 1 & 0\\
                 0 & 1
                \end{array} \right) & 1 & 1 & 1\\
\hline
\mathcal{C}_2 & c_2 = \left ( \begin{array}{rr}
                 1 & 1\\
                 0 & 1
                \end{array} \right) & 1 & \frac{q^2-1}{2} & p\\
\hline
\mathcal{C}_3 & c_3 = \left ( \begin{array}{rr}
                 1 & \Delta\\
                 0 & 1
                \end{array} \right) & 1 & \frac{q^2-1}{2} & p\\
\hline
\mathcal{C}_4 & c_4 = \left ( \begin{array}{rr}
                 x & {\Delta}y\\
                 y & x
                \end{array} \right) x^2-{\Delta}y^2=1,y\neq 0 & \frac{q-1}{4} & q(q-1) & d\mid (q+1)/2,d\neq 2\\
\hline
\mathcal{C}_5 & c_5 = \left ( \begin{array}{rr}
                 i & 0\\
                 0 & i
                \end{array} \right) & 1 & \frac{q(q+1)}{2} & 2\\
\hline
\mathcal{C}_6 & c_6 = \left ( \begin{array}{rr}
                 x & 0\\
                 0 & x^{-1}
                \end{array} \right)  x \not \in \{\pm 1 , \pm i\}  & \frac{q-5}{4} & q(q+1) & d\mid (q-1)/2,d\neq 2\\
\hline
\end{array}
$$
\caption{Conjugacy Classes of $\PSL(2,q)$ for $q \equiv 1\pmod 4$}\label{tpsl1}
\end{table}

\begin{table}[ht]
$$
\begin{array}{|c|c|c|c|c|}
\hline
 & \text{Representative} & \text{No.} & \text{Size} & \text{Order representative} (d)\\
\hline
\mathcal{C}_1 &c_1 = \left ( \begin{array}{rr}
                 1 & 0\\
                 0 & 1
                \end{array} \right) & 1 & 1 & 1\\
\hline
\mathcal{C}_2 & c_2 = \left ( \begin{array}{rr}
                 1 & 1\\
                 0 & 1
                \end{array} \right) & 1 & \frac{q^2-1}{2} & p\\
\hline
\mathcal{C}_3 & c_3 = \left ( \begin{array}{rr}
                 1 & \Delta\\
                 0 & 1
                \end{array} \right) & 1 & \frac{q^2-1}{2} & p\\
\hline
\mathcal{C}_4 & c_4 = \left ( \begin{array}{rr}
                 x & {\Delta}y\\
                 y & x
                \end{array} \right) x^2-{\Delta}y^2=1, xy\neq 0 & \frac{q-3}{4} & q(q-1) & d\mid (q+1)/2,d\neq 2\\
\hline
\mathcal{C}_5 & c_5 = \left ( \begin{array}{rr}
                 0 & -1\\
                 1 & 0
                \end{array} \right) & 1 & \frac{q(q-1)}{2} & 2\\
\hline
\mathcal{C}_6 & c_6 = \left ( \begin{array}{rr}
                 x & 0\\
                 0 & x^{-1}
                \end{array} \right)  x \not \in \{\pm 1 \} & \frac{q-3}{4} & q(q+1) & d\mid (q-1)/2,d\neq 2 \\
\hline
\end{array}
$$
\caption{Conjugacy Classes of $\PSL(2,q)$ for $q \equiv 3\pmod 4$}\label{tpsl2}
\end{table}
From Table \ref{tpsl1} and \ref{tpsl2}, we get that $\PSL(2,q)$ has $(q+5)/2$ conjugacy classes in both cases. We recall the definition of a real conjugacy class below.

\begin{definition}
Let $G$ be a finite group and $C$ be a conjugacy class of $G$. If  $x$ and $x^{-1}$ both are in $C$, then $C$ is called as real conjugacy class. 
\end{definition}

From Lemma \ref{c5s5.6}, we get that the non-involutive real conjugacy classes do not give primitive quandles. The following lemma counts the number of real conjugacy classes in $\PSL(2,q)$.

\begin{lemma}
The group $\PSL(2,q)$ has $(q+5)/2$ real conjugacy classes for $q \equiv 1\pmod 4$ and $(q+1)/2$ real conjugacy classes for $q \equiv 3\pmod 4$.
\end{lemma}
\begin{proof}
The proof follows from \cite[Theorem 4.3]{GS11}.
\end{proof}

The remaining two non real conjugacy classes for $q \equiv 3\pmod 4$ are $\mathcal{C}_2$ and $\mathcal{C}_3$. In fact we have $\mathcal{C}^{-1}_2 = \mathcal{C}_3$. However, in this case, there exist nontrivial $k'$s such that $c_2$ is conjugate to $c^{k}_2$. In fact this happens for all $k\in \mathbb{F}^{2}_q$. As we have
$$ \begin{pmatrix}
              1 &  1\\
              0 &  1
             \end{pmatrix}^k = \begin{pmatrix}
              1 &  k\\
              0 &  1
             \end{pmatrix} \hspace{0.5cm}  \text{and}    \hspace{0.5cm}  \begin{pmatrix}
              a &  b\\
              0 &  a^{-1}
             \end{pmatrix} \begin{pmatrix}
              1 &  1\\
              0 &  1
             \end{pmatrix} \begin{pmatrix}
              a &  b\\
              0 &  a^{-1}
             \end{pmatrix}^{-1} = \begin{pmatrix}
              1 &  a^2\\
              0 &  1
             \end{pmatrix}  $$

Thus, using lemma \ref{c5s5.6}, all conjugacy classes of $\PSL(2,q)$ except those of involution are discarded as candidates for primitive quandles. Now we compute the centralizer of an involution\index{Involution} in $\PSL(2,q)$.

\begin{lemma}
Let $G=\PSL(2,q)$ and $C_G(x)$ be the centralizer\index{Centralizer} of element $x \in G$. Then we have the following:
\begin{enumerate}
\item[$({\rm i})$.]
For $q \equiv 1\pmod 4$, $C_{G}(c_5) \cong D_{(q-1)/2}$
\item[$({\rm ii})$.]
For $q \equiv 3\pmod 4$, $C_{G}(c_5) \cong D_{(q+1)/2}$.

\end{enumerate}
\end{lemma}

\begin{proof}
For part $(\rm i)$, we see that $C_{G}(c_5)$ is generated by $ \left\{r= \begin{pmatrix}
x & 0\\
0 & x^{-1} 
\end{pmatrix},
s=\begin{pmatrix}
0 & 1\\
-1 & 0 
\end{pmatrix} \right\}
$ where $x$ has order $(q-1)/2$ in $\mathbb{F}^{*}_q$. We have $C_{G}(c_5) = \{r,s~|~r^{(q-1)/2},s^2,srs=r^{-1}\}$. By the orbit stabilizer theorem\index{Orbit stabilizer theorem}, we observe that $|C_{G}(c_5)|=q-1$. Thus $C_{G}(c_5)$ is dihedral group of order $q-1$. 

For part $({\rm ii})$, we see that $C_{G}(c_5)$ is generated by $ \left\{r= \begin{pmatrix}
x & -y\\
y & x 
\end{pmatrix},
s=\begin{pmatrix}
a & b\\
b & -a 
\end{pmatrix} \right\}
$ where $x,y \in \mathbb{F}_q $ are chosen such that $x^2+y^2 \equiv 1\pmod p$ and $r$ has order $(q+1)/2$ and $a,b \in \mathbb{F}_q$ such that $a^2+b^2 \equiv 1\pmod p$. We have $C_{G}(c_5) = \{r,s~|~r^{(q+1)/2},s^2,srs=r^{-1}\}$. The orbit stabilizer theorem shows that $|C_{G}(c_5)|=q+1$. Thus $C_{G}(c_5)$ is dihedral group of order $q+1$.

\end{proof}

We note that $\PSL(2,3) \cong A_4$, $\PSL(2,5) \cong A_5$. $C_G(c_5)$ is not maximal subgroup of $\PSL(2,7)$ or $\PSL(2,9)$. The list of maximal subgroups\index{Maximal subgroup} of $\PSL(2,q)$ is available at \cite[Corollary 2.2]{KingBCC05}. We observe from the list that the subgroups $D_{(q-1)/2}$ and $D_{(q+1)/2}$ are maximal subgroups of $\PSL(2,q)$ for $q>9$. Now we state the following theorem.

\begin{theorem}
Let $X$ be primitive quandle with $\operatorname{Inn}(X)=\PSL(2,q)$ for $q>9$ then either:
\begin{itemize}
\item[$({\rm i})$.]
$q \equiv 1\pmod 4, |X|=q(q+1)/2$ and $X \cong \operatorname{Conj}(\PSL(2,q),\mathcal{C}_5)$ or
\item[$({\rm ii})$.]
$q \equiv 3\pmod 4, |X|=q(q-1)/2$ and $X \cong \operatorname{Conj}(\PSL(2,q),\mathcal{C}_5)$.

\end{itemize}
\end{theorem}

\begin{proof}
In both case, the conjugacy class $\mathcal{C}_5$ forms conjugation quandle $X=\operatorname{Conj}(\PSL(2,q),C_5)$ with $\operatorname{Inn}(X)=\PSL(2,q)$. The stabilizer\index{Stabilizer} of an $c_5 \in C_5$ under the conjugation action is the centralizer, which is a maximal subgroup. Hence, the action of $\operatorname{Inn}(X)$ on $X$ is primitive, so the quandle is primitive.
\end{proof}

Now, we describe primitive quandles arising from $\PGL(2,q)$. The conjugacy classes of $\PGL(2,q)$ are given in Table \ref{tpgl}. [\cite{Ada02},\cite{FGV10}]

\begin{table}[ht]
$$
\begin{array}{|c|c|c|c|c|}
\hline
 & \text{Representative} & \text{No.} & \text{Size} & \text{Order representative} (d)\\
\hline
\mathcal{C}_1' &c_1' = \left ( \begin{array}{rr}
                 1 & 0\\
                 0 & 1
                \end{array} \right) & 1 & 1 & 1\\
\hline
\mathcal{C}_2' & c_2' = \left ( \begin{array}{rr}
                 1 & 1\\
                 0 & 1
                \end{array} \right) & 1 & q^2-1 & p\\
\hline
\mathcal{C}_3' & \quad c_3' =  \left ( \begin{array}{rr}
                 x & 0\\
                 0 & 1
                \end{array} \right),  x \neq \pm 1  \quad \quad & \frac{q-3}{2} & ~~ q(q+1) & d\mid q-1,d\neq 2\\
\hline
\mathcal{C}_4' & c_4' = \left ( \begin{array}{rr}
                 -1 & 0\\
                 0 & 1
                \end{array} \right) & 1 & \frac{q(q+1)}{2} & 2\\
\hline
\mathcal{C}_5' & c_5' = \left ( \begin{array}{rr}
                 x & {\Delta}y\\
                 y & x
                \end{array} \right),  xy \neq 0 \quad & \frac{q-1}{2} & q(q-1) & d\mid q+1,d\neq 2\\
\hline
\mathcal{C}_6' & c_6' = \left ( \begin{array}{rr}
                 0 & \Delta\\
                 1 & 0
                \end{array} \right) & 1 & \frac{q(q-1)}{2} & 2\\
\hline
\end{array}
$$
\caption{Conjugacy Classes of $\PGL(2,q)$}\label{tpgl}
\end{table}

Table \ref{tpgl} shows that $\PGL(2,q)$ has $q+2$ conjugacy classes. Using Theorem \cite[Theorem 2.8]{GS11}, we get that all conjugacy classes of $\PGL(2,q)$ are real. Thus, using lemma \ref{c5s5.6}, all conjugacy classes of $\PGL(2,q)$ except those of involution are discarded as candidates for primitive quandles.

\begin{lemma}
Let $x \in \PGL(2,q)$ be such that $\operatorname{det}(x)$ is a square element. Then $x\in \PSL(2,q)$.
\end{lemma}

We get that $\mathcal{C}_4'$ and $\mathcal{C}_6'$ are conjugacy classes of order two elements. As a consequence of above lemma, $q \equiv 1\pmod 4$, $\mathcal{C}_4'$ is in $\PSL(2,q)$ and for $q \equiv 3\pmod 4$, $\mathcal{C}_6'$ is in $\PSL(2,q)$. Hence we have $\langle \mathcal{C}_6' \rangle = \PGL(2,q)$ for $q \equiv 1\pmod 4$ and $\langle \mathcal{C}_4' \rangle = \PGL(2,q)$ for $q \equiv 3\pmod 4$. Now we compute the centralizer of an involution in $\PGL(2,q)$.

\begin{lemma}
Let $G=\PGL(2,q)$ and $C_G(x)$ be the centralizer of $x \in G$. Then we have the following:
\begin{enumerate}
\item[$({\rm i})$.]
$C_{G}(c_4') \cong D_{q-1}$
\item[$({\rm ii})$.]
$C_{G}(c_6') \cong D_{q+1}$.
\end{enumerate}
\end{lemma}

\begin{proof}
For part $({\rm i})$, we get that the centralizer of involution $c_4$ is dihedral group of order $2q-2$ generated by $\left\{ r= \begin{pmatrix}
x & 0\\
0 & 1 
\end{pmatrix}, s= \begin{pmatrix}
0 & 1\\
1 & 0 
\end{pmatrix}    \right\} $ where $x \in \mathbb{F}^{*}_q$ is such that $r$ has maximum order $q-1$. For part $({\rm ii})$, we get that the centralizer of involution $c_6$ is dihedral group of order $2q+2$ generated by $\left\{ r= \begin{pmatrix}
x & {\Delta}y\\
y & x 
\end{pmatrix}, s= \begin{pmatrix}
1 & 0\\
0 & -1 
\end{pmatrix}    \right\} $ where $r$ has maximum order $q+1$ which is possible due to embedding of $F^{*}_{q^2}$ in $\GL(2,q)$. 

\end{proof}

We have $\PGL(2,3) \cong S_4$. The list of maximal subgroups of $\PGL(2,q)$ is available at \cite[Corollary 2.3]{KingBCC05}. We observe from the list that the subgroups $D_{q-1}$ and $D_{q+1}$ are maximal subgroups of $\PGL(2,q)$ for $q>3$. Thus, we get the following theorem.

\begin{theorem}
Let $X$ be primitive quandle with $\operatorname{Inn}(X)=\PGL(2,q)$ with $q>3$ then either:
\begin{enumerate}
\item[$({\rm i})$.]
$q \equiv 1\pmod 4, |X|=q(q-1)/2$ and $X \cong \operatorname{Conj}(\PGL(2,q),\mathcal{C}_6')$ or
\item[$({\rm ii})$.]
$q \equiv 3\pmod 4, |X|=q(q+1)/2$ and $X \cong \operatorname{Conj}(\PGL(2,q),\mathcal{C}_4')$.

\end{enumerate}

\end{theorem}

\begin{proof}
We note that $\mathcal{C}_6'$ and $\mathcal{C}_4'$ generate $\PGL(2,q)$ in case $({\rm i})$ and $({\rm ii})$ respectively. Moreover, the action of $\operatorname{Inn}(X)$ on $X$ is a primitive action as the centralizers are maximal subgroups. Thus, the quandles are primitive in both cases.
\end{proof}

In even characteristic, we have $\PGL(2,2^n) \cong \PSL(2,2^n) \cong \SL(2,2^n)$. The conjugacy classes of $\SL(2,2^n)$ are given below.

\begin{table}[ht]
$$
\begin{array}{|c|c|c|c|c|}
\hline
 & \text{Representative} & \text{No.} & \text{Size} & \text{Centralizer}\\
\hline
\mathcal{C}_1 &c_1 = \left ( \begin{array}{rr}
                 1 & 0\\
                 0 & 1
                \end{array} \right) & 1 & 1 & SL(2,2^n)\\
\hline
\mathcal{C}_2 & c_2 = \left ( \begin{array}{rr}
                 1 & 1\\
                 0 & 1
                \end{array} \right) & 1 & 2^{2n}-1 & \mathbb{F}_{2^n}\\
\hline
\mathcal{C}_3 & \quad c_3 =  \left ( \begin{array}{rr}
                 x & 0\\
                 0 & x^{-1}
                \end{array} \right) \left( x \neq 1 \right) \quad \quad & \frac{2^n-2}{2} & ~~ 2^n(2^n+1) & \mathbb{Z}/(2^n-1)\mathbb{Z}\\
\hline
\mathcal{C}_4 & c_4 = \left ( \begin{array}{rr}
                 0 & 1\\
                 1 & x+\bar{x}
                \end{array} \right) x \in \mathbb{F}_{2^{2n}}\setminus \mathbb F_{2^n} & \frac{2^n}{2} & 2^n(2^n-1) & \mathbb{Z}/(2^n+1)\mathbb{Z}\\
\hline
\end{array}
$$
\caption{Conjugacy Classes of $\SL(2,2^n)$}\label{tsl}
\end{table}

We observe using \cite[Corollary 2.2]{KingBCC05} that centralizers are not maximal subgroups in $\SL(2,2^n)$, so no conjugacy class gives a primitive quandle. Now we describe the isomorphism classes and profiles of simple quandles arising from $\PSL(2,q)$ and $\PGL(2,q)$. 

\subsection{Isomorphism classes} We describe the orbits of the action of ${\rm Out}(G)$ on the set of conjugacy classes. Using Theorem \ref{isoq}, this gives us an enumeration of simple quandles up to isomorphism. We note that ${\rm Out}(PSL(2,q))= \langle \phi,\rho \rangle \cong \mathbb{Z}_2 \times \mathbb{Z}_n$, where $\phi(x)=DxD^{-1}$ for $D=\left ( \begin{array}{rr}
                 \Delta & 0\\
                 0 & 1
                \end{array}\right) $
and $\rho \left( \begin{array}{rr}
                 a & b\\
                 c & d
                \end{array}\right)=\left( \begin{array}{rr}
                 a^p & b^p\\
                 c^p & d^p
                \end{array}\right).$ We get $\phi(\mathcal C_2)=\mathcal{C}_3$. Further, $\phi$ maps classes of type $\mathcal{C}_4,\mathcal{C}_5$ and $\mathcal{C}_6$ to themselves.
In $\mathcal{C}_4$ (or $\mathcal{C}_6$), for each divisor $d>2$ of $\frac{q+1}{2}$ (or $\frac{q-1}{2}$), there are $\varphi(d)/2$ classes. For each $d>2$, the size of the orbit of a corresponding class under $\langle \rho \rangle$ is the minimum $r>0$ such that $p^r \equiv \pm 1 \pmod d$.

For $G=\PGL(2,q)$, we have ${\rm Out}(G) = \langle \rho \rangle \cong \mathbb{Z}_n$. Elements of ${\rm Out}(G)$ map $\mathcal C_2',\mathcal C_4',\mathcal C_6'$ to themselves. For $\mathcal{C}_3'$ (or $\mathcal{C}_5'$), let $d$ be the divisor of $q-1$ (or $q+1$). If $d \nmid \frac{q-1}{2}$ (or $d \nmid \frac{q+1}{2}$) then the corresponding class generates $\PGL(2,q)$. The sizes of orbits is the same as in the above case.

For $\SL(2,2^n)$, we have ${\rm Out}(SL(2,2^n)) = \langle \rho \rangle \cong \mathbb{Z}/n\mathbb{Z}$ consisting of field automorphisms. The sizes of orbits of the action of $\langle \rho \rangle$ on the set of Conjugacy classes is similar to that done in $\PSL(2,q)$.

\subsection{Description of \textit{profile}} We observe from the data of profiles of simple quandles ${\rm Conj}(G,C)$ for $G=\PSL(2,q)$ and $\PGL(2,q)$ that cycle structures have only two distinct lengths. We prove this in general below.

\begin{lemma}\label{prof}
Let $G=\PSL(2,q)$ or  $\PGL(2,q), q$ odd, and $C$ be a nontrivial conjugacy class of $G$. Let $x,y \in C$ then either $x$ commutes with $y$ or $\langle x \rangle \cap C_G(y)=\{1\}$.
\end{lemma}

\begin{proof}
We begin with $G=\PSL(2,q)$. Let $x\in \mathcal C_2,\mathcal C_3$ or $\mathcal C_5$, then $o(x)$ is prime, so for any $y$ in the same conjugacy class, if $x^m\neq 1$ commutes with $y$, then $\langle x \rangle \subseteq C_G(y)$. Hence, the result holds. For conjugacy class of type $C_4$, we have $C_{G}(c_4)=\left\{  \begin{pmatrix}
a & \Delta b\\
b & a
\end{pmatrix} \;|\; a^2-\Delta b^2=1 \right\} /\{\pm I\} \cong \mathbb{Z}/\frac{(q+1)}{2}\mathbb Z$. For conjugacy class of type $\mathcal C_6$, we have 
$C_{G}(c_6)=\left\{  \begin{pmatrix}
                 a & 0\\
                 0 & a^{-1}
                \end{pmatrix} \;|\; a \in \mathbb{F}_q^{*} \right\} /\{\pm I\} \cong \mathbb{Z}/\frac{(q-1)}{2}\mathbb Z$.

Let $q\equiv 1 \pmod 4$,  $x,y \in \mathcal C_6$ and $g \in C_G(x) \cap C_G(y)$. As $C_G(x)$ and $C_G(y) $ are abelian, we get $C_G(x),C_G(y)\leq C_G(g)$.  If $o(g)>2$, then $g$ is of type $\mathcal C_6$ and so $C_G(x)=C_G(y)$. If $o(g)=2$ then $C_G(g) \cong D_{(q- 1)/2}$. Since $C_G(g)$ contains a unique cyclic subgroup of size $(q-1)/2$, so $C_G(x)=C_G(y)$. The $o(g)=2$ is not possible for $\mathcal C_4$ as $C_G(x)$ has odd order. The case of $q \equiv 3 \pmod 4$ is similar.

The case of $G=\PGL(2,q)$ is the same as above. We have $C_{G}(c_3')=\left\{  \begin{pmatrix}
                 a & 0\\
                 0 & 1
                \end{pmatrix} \;|\; a\in F_{q}^* \right\} \cong \mathbb{Z}/(q-1)\mathbb Z$ and $C_{G}(c_5')=\left\{  \begin{pmatrix}
                 a & \Delta b\\
                 b & a
                \end{pmatrix} \;|\; a^2-\Delta b^2\neq 0 \right\} /\{\pm I\} \cong \mathbb{Z}/(q+1)\mathbb Z$. If $x,y \in \mathcal C_3'$ ( or $\mathcal C_5'$), then we have either $C_G(x)\cap C_G(y)=\{1\}$ or $C_G(x)=C_G(y)$. Hence, the result follows.
\end{proof}

\begin{theorem}
Let ${\rm Conj}(G,C)$ be simple quandle arising from $G=\PSL(2,q), \PGL(2,q), q$ odd, then the \textit{profile} of ${\rm Conj}(G,C)$ is of form $1^k\lambda^l$ where $\lambda={\rm ord}(x)$ for $x\in C$.
\end{theorem}

\begin{proof}
Let $x \in C$ and $R_x$ be the right translation map. Then for $y\in C$, the cycle in $R_x$ containing $y$ is $(y,xyx^{-1},\dots,x^{t}yx^{-t})$. The result now follows from Lemma \ref{prof}.
\end{proof}

A finite group is called a CA-group if the centralizer of every noncentral element is abelian. We note that ${\SL}(2,2^n)$ is CA-group.

\begin{theorem}
Let ${\rm Conj}(G,C)$ be simple quandle for $G={\SL}(2,2^n)$ then the \textit{profile} of ${\rm Conj}(G,C)$ is of form $1^k\lambda^l$.
\end{theorem}
\begin{proof}
Let $x,y \in C$, since $G$ is a CA-group so either $C_G(x)=C_G(y)$ or $C_G(x) \cap C_G(y)=\{1\}$ from \cite[Lemma 2.1]{Bro18}. Hence, the \textit{profile} is of the form $1^k\lambda^l$ where $\lambda={\rm ord}(x)$.
\end{proof}

We observe from the data that even after excluding the above cases, many simple quandles have \textit{profile} $1^k\lambda^l$. Therefore, we state the following problem.

\begin{problem}
Classify finite simple quandles with \textit{profile} $1^k\lambda^l$.
\end{problem}

\newpage
\begin{table}[htb]
\centering
\caption{ Non primitive simple quandles up to order 4096}\label{tables1}
\begin{small}
\begin{tabular}{@{}llll@{\hskip 4mm}|@{\hskip 4mm}llll@{}}
\toprule
$X$ &$\operatorname{Inn}(X)$ &$\operatorname{Dis}(X)$ & $\rm{Prof}(X)$         &$X$  &$\operatorname{Inn}(X)$  &$\operatorname{Dis}(X)$ & $\rm{Prof}(X)$\\
\midrule
$[12,1]$ & $A_5$ & $A_5$                                               & $1^2 5^2$                                                                              &$[156,4]$ & $\PSL(2,13) \rtimes \mathbb{Z}_2$ & $\PSL(2,13)$   & $1^2 14^{11}$                                 \\        
 $[15,1]$ & $A_5$ & $A_5$                                              & $1^3 2^6$                                                                     &$[156,5]$ & $\PSL(2,13) \rtimes \mathbb{Z}_2$ & $\PSL(2,13)$            & $1^2 14^{11}$                              \\         
 $[20,1]$ & $A_5$ & $A_5$                                              & $1^2 3^6$                                                               &$[156,6]$ & $\PSL(2,13) \rtimes \mathbb{Z}_2$ & $\PSL(2,13)$                  & $1^2 14^{11}$                              \\         
 $[20,2]$ & $S_5$ & $A_5$                                              & $1^2 6^3$                                                         &$[168,1]$ & $A_{9}$ & $A_{9}$                                                     & $1^{42} 3^{42}$                              \\        
 $[21,1]$ & $\PSL(3,2)$ & $\PSL(3,2)$                                    & $1^5 2^8$                                                  &$[168,2]$ & $\PSL(2,8) \rtimes \mathbb{Z}_3$ & $\PSL(2,8)$                                 & $1^3 3^1 9^{18}$                              \\         
 $[24,1]$ & $\PSL(3,2)$ & $\PSL(3,2)$                                    & $1^3 7^3$                                            &$[168,3]$ & $\PSL(2,8)\rtimes \mathbb{Z}_3$ & $\PSL(2,8)$                                        & $1^3 3^1 9^{18}$                              \\
\bottomrule

\end{tabular}
\end{small}
\end{table}
\newpage

\begin{table}[htb]
\centering
\caption{Continued}\label{tables2}
\begin{small}
\begin{tabular}{@{}llll@{\hskip 4mm}|@{\hskip 4mm}llll@{}}
\toprule
$X$ &$\operatorname{Inn}(X)$ &$\operatorname{Dis}(X)$ & $\rm{Prof}(X)$         &$X$  &$\operatorname{Inn}(X)$  &$\operatorname{Dis}(X)$ & $\rm{Prof}(X)$\\
\midrule
 $[312,2]$ & $O(5,5)$ & $O(5,5)$                                  & $1^{62} 5^{50}$                                                                                    &$[506,7]$ & $\PSL(2,23) \rtimes \mathbb{Z}_2$ & $\PSL(2,23)$                        & $1^2 24^{21}$                          \\        
 $[315,1]$ & $\PSL(3,4)$ & $\PSL(3,4)$                              & $1^{27} 2^{144}$                                                                                &$[506,8]$ & $\PSL(2,23) \rtimes \mathbb{Z}_2$ & $\PSL(2,23)$                           & $1^2 8^{63}$                          \\         
 $[330,1]$ & $A_{11}$ & $A_{11}$                                  & $1^{114} 3^{72}$                                                                                   &$[506,9]$ & $\PSL(2,23) \rtimes \mathbb{Z}_2$ & $\PSL(2,23)$                          & $1^2 24^{21}$                          \\         
 $[336,1]$ & $\PSL(3,4)\rtimes \mathbb{Z}_3$ & $\PSL(3,4)$          & $1^{21} 3^{105}$                                                                 &$[506,10]$ & $\PSL(2,23) \rtimes \mathbb{Z}_2$ & $\PSL(2,23)$                                         & $1^2 24^{21}$                        \\        
 $[342,1]$ & $\PSL(2,19)$ & $\PSL(2,19)$                            & $1^2 5^{68}$                                                                              &$[506,11]$ & $\PSL(2,23) \rtimes \mathbb{Z}_2$ & $\PSL(2,23)$                                & $1^2 24^{21}$                             \\
\bottomrule
\end{tabular}
\end{small}
\end{table}
\newpage

\begin{table}[htb]
\centering
\caption{Continued}\label{tables3}
\begin{small}
\begin{tabular}{@{}llll@{\hskip 4mm}|@{\hskip 4mm}llll@{}}
\toprule
$X$ &$\operatorname{Inn}(X)$ &$\operatorname{Dis}(X)$ & $\rm{Prof}(X)$         &$X$  &$\operatorname{Inn}(X)$  &$\operatorname{Dis}(X)$ & $\rm{Prof}(X)$\\
\midrule
 $[720,1]$ & $O(5,3)$ & $O(5,3)$                                           & $1^4 2^1 3^{20} 6^{109}$                                                                           &$[930,1]$ & $\PSL(2,31) $ & $\PSL(2,31)$                         & $1^2 16^{58}$                  \\        
 $[720,2]$ & $(A_5 \times A_5) \rtimes \mathbb{Z}_2$ & $A_5 \times A_5$    & $1^2 2^4 5^2 10^{70}$                                          &$[930,2]$ & $\PSL(2,31) $ & $\PSL(2,31)$                                                             & $1^2 8^{116}$             \\         
 $[720,3]$ & $M_{11}$ & $M_{11}$                                           & $1^5 11^{65}$                                                                        &$[930,3]$ & $\PSL(2,31) $ & $\PSL(2,31)$                                       & $1^2 16^{58}$                      \\         
 $[720,4]$ & $M_{11}$ & $M_{11}$                                           & $1^5 11^{65}$                                                                        &$[930,4]$ & $\PSL(2,31) $ & $\PSL(2,31)$                                       & $1^2 4^{232}$                   \\        
 $[728,1]$ & $A_{14}$ & $A_{14}$                                           & $1^{332} 3^{132}$                                                                         &$[930,5]$ & $\PSL(2,31) $ & $\PSL(2,31)$                                  & $1^2 16^{58}$                      \\         
 $[728,2]$ & $G_{2}(3)$ & $G_{2}(3)$                                       & $1^{26} 3^{234}$                                                                       &$[930,6]$ & $\PSL(2,31)$ & $\PSL(2,31)$                                      & $1^2 16^{58}$                        \\
\bottomrule
\end{tabular}
\end{small}
\end{table}
\newpage

\begin{table}[htb]
\centering
\caption{Continued}\label{tables4}
\begin{small}
\begin{tabular}{@{}llll@{\hskip 4mm}|@{\hskip 4mm}llll@{}}
\toprule
$X$ &$\operatorname{Inn}(X)$ &$\operatorname{Dis}(X)$ & $\rm{Prof}(X)$         &$X$  &$\operatorname{Inn}(X)$  &$\operatorname{Dis}(X)$ & $\rm{Prof}(X)$\\
\midrule
  $[1056,1]$ & $\PSL(2,32)$ & $\PSL(2,32)$                                     & $1^{2} 31^{34}$                                                                            &$[1406,1]$ & $\PSL(2,37) $ & $\PSL(2,37)$                                                  &$1^2 9^{156}$                                           \\        
 $[1056,2]$ & $\PSL(2,32)$ & $\PSL(2,32)$                                     & $1^2 31^{34}$                                                                            &$[1406,2]$ & $\PSL(2,37) $ & $\PSL(2,37)$                                                    &$1^2 18^{78}$                                      \\         
 $[1056,3]$ & $\PSL(2,32)$ & $\PSL(2,32)$                                     & $1^2 31^{34}$                                                                            &$[1406,3]$ & $\PSL(2,37) $ & $\PSL(2,37)$                                                    &$1^2 18^{78}$                                    \\         
 $[1080,1]$ & $\PSL(4,3) \rtimes \mathbb{Z}_2$ & $\PSL(4,3)$                  & $1^{118} 2^{281}$                                                              &$[1406,4]$ & $\PSL(2,37) $ & $\PSL(2,37)$                                                              &$1^2 9^{156}$                                      \\        
 $[1104,1]$ & $\PSL(2,47)$ & $\PSL(2,47)$                                     & $1^{23} 47^{23}$                                                                            &$[1406,5]$ & $\PSL(2,37) $ & $\PSL(2,37)$                                                 &$1^2 18^{78}$                                         \\
\bottomrule
\end{tabular}
\end{small}
\end{table}

\newpage

\begin{table}[htb]
\centering
\caption{Continued}\label{tables5}
\begin{small}
\begin{tabular}{@{}llll@{\hskip 4mm}|@{\hskip 4mm}llll@{}}
\toprule
$X$ &$\operatorname{Inn}(X)$ &$\operatorname{Dis}(X)$ & $\rm{Prof}(X)$         &$X$  &$\operatorname{Inn}(X)$  &$\operatorname{Dis}(X)$ & $\rm{Prof}(X)$\\
\midrule
$[1722,4]$ & $\PSL(2,41)$ & $\PSL(2,41)$                                       &$1^2 4^{430}$                                                                          &$[1892,7]$ & $\PSL(2,43) $ & $\PSL(2,43)$                                                   &$1^2 7^{270}$                   \\        
$[1722,5]$ & $\PSL(2,41)$ & $\PSL(2,41)$                                       &$1^2 20^{86}$                                                                          &$[1892,8]$ & $\PSL(2,43) $ & $\PSL(2,43)$                                                   &$1^2 21^{90}$              \\         
$[1722,6]$ & $\PSL(2,41)$ & $\PSL(2,41)$                                       &$1^2 20^{86}$                                                                          &$[1892,9]$ & $\PSL(2,43) $ & $\PSL(2,43)$                                                   &$1^2 21^{90}$               \\         
$[1722,7]$ & $\PSL(2,41)$ & $\PSL(2,41)$                                       &$1^2 5^{344}$                                                                          &$[1892,10]$ & $\PSL(2,43) $ & $\PSL(2,43)$                                                  &$1^2 21^{90}$                       \\        
$[1722,8]$ & $\PSL(2,41)$ & $\PSL(2,41)$                                       &$1^2 10^{172}$                                                                          &$[1892,11]$ & $\PSL(2,43)\rtimes \mathbb{Z}_2 $ & $\PSL(2,43)$                             &$1^2 42^{45}$                       \\
\bottomrule
\end{tabular}
\end{small}
\end{table}

\newpage

\begin{table}[htb]
\centering
\caption{Continued}\label{tables6}
\begin{small}
\begin{tabular}{@{}llll@{\hskip 4mm}|@{\hskip 4mm}llll@{}}
\toprule
$X$ &$\operatorname{Inn}(X)$ &$\operatorname{Dis}(X)$ & $\rm{Prof}(X)$         &$X$  &$\operatorname{Inn}(X)$  &$\operatorname{Dis}(X)$ & $\rm{Prof}(X)$\\
\midrule
$[2240,1]$ & $Sz(8)$ & $Sz(8)$                                    &$1^4 13^{172}$                                                                             &$[2450,9]$ &  $\PSL(2,49)\rtimes \mathbb{Z}_2$ & $\PSL(2,49)$                       &$1^2 48^{51}$                                              \\        
$[2240,2]$ & $A_9$ & $A_9$                                        &$1^{26} 3^{738}$                                                                             &$[2450,10]$ & $\PSL(2,49)\rtimes \mathbb{Z}_2 $ & $\PSL(2,49)$                    &$1^2 48^{51}$                                              \\         
$[2240,3]$ & $O^{+}(8,2)$ & $O^{+}(8,2)$                          &$1^{242} 3^{666}$                                                                      &$[2450,11]$ & $\PSL(2,49)\rtimes \mathbb{Z}_2 $ & $\PSL(2,49)$                          &$1^2 48^{51}$                                               \\         
$[2240,4]$ & $\PSL(3,4)$ & $\PSL(3,4)$                              &$1^8 3^{744}$                                                                      &$[2450,12]$ & $\PSL(2,49)\rtimes \mathbb{Z}_2 $ & $\PSL(2,49)$                              &$1^2 16^{153}$                                                \\        
$[2240,5]$ & $\PSL(3,4) \rtimes \mathbb{Z}_2$ & $\PSL(3,4)$         &$1^8 6^{372}$                                                      &$[2450,13]$ & $\PSL(2,49)\rtimes \mathbb{Z}_2 $ & $\PSL(2,49)$                                              &$1^2 16^{153}$                                                 \\
\bottomrule
\end{tabular}
\end{small}
\end{table}
\newpage
\begin{table}[htb]
\centering
\caption{Continued}\label{tables7}
\begin{small}
\begin{tabular}{@{}llll@{\hskip 4mm}|@{\hskip 4mm}llll@{}}
\toprule
$X$ &$\operatorname{Inn}(X)$ &$\operatorname{Dis}(X)$ & $\rm{Prof}(X)$         &$X$  &$\operatorname{Inn}(X)$  &$\operatorname{Dis}(X)$ & $\rm{Prof}(X)$\\
\midrule
 $[2793,3]$ & $\PSL(3,7) \rtimes \mathbb{Z}_3$ & $\PSL(3,7)$           &$1^{57} 6^{456}$                                              &$[3360,1]$   &  $A_8$ & $A_8$                                              &$1^2 2^2 3^{10} 6^{554}$                                                            \\        
 $[2862,1]$ & $\PSL(2,53)$ & $\PSL(2,53)$                              &$1^2 26^{110}$                                                &$[3360,2]$   &  $S_8$ & $A_8$                                              &$1^2 2^2 3^2 6^{558}$                                                  \\         
 $[2862,2]$ & $\PSL(2,53)$ & $\PSL(2,53)$                              &$1^2 13^{220}$                                                &$[3360,3]$   &  $S_8$ & $A_8$                                              &$1^4 3^4 4^{14} 12^{274}$                                                   \\         
 $[2862,3]$ & $\PSL(2,53)$ & $\PSL(2,53)$                              &$1^2 13^{220}$                                                &$[3360,4]$   &  $A_9$ & $A_9$                                              &$1^{12} 3^{1116}$                                                   \\        
$[2862,4]$ & $\PSL(2,53)$ & $\PSL(2,53)$                               &$1^2 13^{220}$                                                &$[3360,5]$   &  $\PSL(3,4)\rtimes \mathbb{Z}_6 $ & $\PSL(3,4)$               &$1^6 2^{12} 3^{14} 6^{548}$                                                      \\
\bottomrule
\end{tabular}
\end{small}
\end{table}

\newpage

\begin{table}[!htb]
\centering
\caption{Continued}\label{tables8}
\begin{small}
\begin{tabular}{@{}llll@{\hskip 4mm}|@{\hskip 4mm}llll@{}}
\toprule
$X$ &$\operatorname{Inn}(X)$ &$\operatorname{Dis}(X)$ & $\rm{Prof}(X)$         &$X$  &$\operatorname{Inn}(X)$  &$\operatorname{Dis}(X)$ & $\rm{Prof}(X)$\\
\midrule
$[3540,13]$ & $\PSL(2,59)$ & $\PSL(2,59)$                                 &$1^2 29^{122}$                                                                    &$[3782,1]$   &  $\PSL(2,61)$ & $\PSL(2,61)$                                    &$1^2 15^{252}$                                         \\        
$[3540,14]$ & $\PSL(2,59)$ & $\PSL(2,59)$                                 &$1^2 29^{122}$                                                                    &$[3782,2]$   &  $\PSL(2,61)$ & $\PSL(2,61)$                                    &$1^2 30^{126}$                                       \\         
$[3540,15]$ & $\PSL(2,59)\rtimes \mathbb{Z}_2$ & $\PSL(2,59)$             &$1^2 58^{61}$                                                                     &$[3782,3]$   &  $\PSL(2,61)$ & $\PSL(2,61)$                                    &$1^2 30^{126}$                                         \\         
$[3540,16]$ & $\PSL(2,59)\rtimes \mathbb{Z}_2$ & $\PSL(2,59)$             &$1^2 58^{61}$                                                                     &$[3782,4]$   &  $\PSL(2,61)$ & $\PSL(2,61)$                                    &$1^2 10^{378}$                                         \\        
$[3540,17]$ & $\PSL(2,59)\rtimes \mathbb{Z}_2$ & $\PSL(2,59)$             &$1^2 58^{61}$                                                                     &$[3782,5]$   &  $\PSL(2,61)$ & $\PSL(2,61)$                                    &$1^2 5^{756}$                                         \\
\bottomrule
\end{tabular}
\end{small}
\end{table}

\newpage

\medskip \noindent \textbf{Data availability.}
The link to the data generated and analyzed is provided in the article.

\medskip \noindent \textbf{Acknowledgement.}
The first-named author would like to acknowledge the support of SERB through the MATRICS project MTR/2022/000231. 

\bigskip

\printbibliography

@article {AG03,
    AUTHOR = {Andruskiewitsch, Nicol\'{a}s and Gra\~{n}a, Mat\'{\i}as},
     TITLE = {From racks to pointed {H}opf algebras},
   JOURNAL = {Adv. Math.},
  FJOURNAL = {Advances in Mathematics},
    VOLUME = {178},
      YEAR = {2003},
    NUMBER = {2},
     PAGES = {177--243},
      ISSN = {0001-8708,1090-2082},
   MRCLASS = {16W30 (17B37 57M27)},
  MRNUMBER = {1994219},
MRREVIEWER = {Ian\ M.\ Musson},
       DOI = {10.1016/S0001-8708(02)00071-3},
       URL = {https://doi.org/10.1016/S0001-8708(02)00071-3},
}

@article {Bor17,
    AUTHOR = {Bors, Alexander},
     TITLE = {Finite groups with an automorphism of large order},
   JOURNAL = {J. Group Theory},
  FJOURNAL = {Journal of Group Theory},
    VOLUME = {20},
      YEAR = {2017},
    NUMBER = {4},
     PAGES = {681--717},
      ISSN = {1433-5883,1435-4446},
   MRCLASS = {20D45},
  MRNUMBER = {3667115},
MRREVIEWER = {Gary\ L.\ Peterson},
       DOI = {10.1515/jgth-2017-0003},
       URL = {https://doi.org/10.1515/jgth-2017-0003},
}

@article {Bro18,
    AUTHOR = {Brough, Julian},
     TITLE = {Central intersections of element centralizers},
   JOURNAL = {Asian-Eur. J. Math.},
  FJOURNAL = {Asian-European Journal of Mathematics},
    VOLUME = {11},
      YEAR = {2018},
    NUMBER = {5},
     PAGES = {1850095, 12},
      ISSN = {1793-5571,1793-7183},
   MRCLASS = {20E34 (20D99)},
  MRNUMBER = {3869872},
MRREVIEWER = {Jiangtao\ Shi},
       DOI = {10.1142/S179355711850095X},
       URL = {https://doi.org/10.1142/S179355711850095X},
}

@article {BS21,
    AUTHOR = {Bonatto, Marco and Stanovsk\'{y}, David},
     TITLE = {Commutator theory for racks and quandles},
   JOURNAL = {J. Math. Soc. Japan},
  FJOURNAL = {Journal of the Mathematical Society of Japan},
    VOLUME = {73},
      YEAR = {2021},
    NUMBER = {1},
     PAGES = {41--75},
      ISSN = {0025-5645,1881-1167},
   MRCLASS = {57K12 (08A30 20N02 20N05)},
  MRNUMBER = {4203321},
MRREVIEWER = {Emanuele\ Zappala},
       DOI = {10.2969/jmsj/83168316},
       URL = {https://doi.org/10.2969/jmsj/83168316},
}

@article {CS25,
    AUTHOR = {Cvr\v{c}ek, Milan and Stanovsk\'{y}, David},
     TITLE = {Primitive quandles with alternating displacement group},
   JOURNAL = {Hiroshima Math. J.},
  FJOURNAL = {Hiroshima Mathematical Journal},
    VOLUME = {55},
      YEAR = {2025},
    NUMBER = {3},
     PAGES = {349--359},
      ISSN = {0018-2079,2758-9641},
   MRCLASS = {20N02 (20B15 57K12)},
  MRNUMBER = {4998225},
       DOI = {10.32917/h2024010},
       URL = {https://doi.org/10.32917/h2024010},
}

@article {GS11,
    AUTHOR = {Gill, Nick and Singh, Anupam},
     TITLE = {Real and strongly real classes in {${\rm PGL}_n(q)$} and
              quasi-simple covers of {${\rm PSL}_n(q)$}},
   JOURNAL = {J. Group Theory},
  FJOURNAL = {Journal of Group Theory},
    VOLUME = {14},
      YEAR = {2011},
    NUMBER = {3},
     PAGES = {461--489},
      ISSN = {1433-5883,1435-4446},
   MRCLASS = {20G40 (20E45)},
  MRNUMBER = {2794378},
MRREVIEWER = {Antonio\ Cossidente},
       DOI = {10.1515/JGT.2010.055},
       URL = {https://doi.org/10.1515/JGT.2010.055},
}

@article {HSV16,
    AUTHOR = {Hulpke, Alexander and Stanovsk\'{y}, David and
              Vojt\v{e}chovsk\'{y}, Petr},
     TITLE = {Connected quandles and transitive groups},
   JOURNAL = {J. Pure Appl. Algebra},
  FJOURNAL = {Journal of Pure and Applied Algebra},
    VOLUME = {220},
      YEAR = {2016},
    NUMBER = {2},
     PAGES = {735--758},
      ISSN = {0022-4049,1873-1376},
   MRCLASS = {57M27 (20B10 20N02)},
  MRNUMBER = {3399387},
MRREVIEWER = {Tetsuya\ Ito},
       DOI = {10.1016/j.jpaa.2015.07.014},
       URL = {https://doi.org/10.1016/j.jpaa.2015.07.014},
}

@article {Joy82a,
    AUTHOR = {Joyce, David},
     TITLE = {A classifying invariant of knots, the knot quandle},
   JOURNAL = {J. Pure Appl. Algebra},
  FJOURNAL = {Journal of Pure and Applied Algebra},
    VOLUME = {23},
      YEAR = {1982},
    NUMBER = {1},
     PAGES = {37--65},
      ISSN = {0022-4049,1873-1376},
   MRCLASS = {57M25 (20F29 20N05 53C35)},
  MRNUMBER = {638121},
MRREVIEWER = {Mark\ E.\ Kidwell},
       DOI = {10.1016/0022-4049(82)90077-9},
       URL = {https://doi.org/10.1016/0022-4049(82)90077-9},
}

@article {Joy82b,
    AUTHOR = {Joyce, David},
     TITLE = {Simple quandles},
   JOURNAL = {J. Algebra},
  FJOURNAL = {Journal of Algebra},
    VOLUME = {79},
      YEAR = {1982},
    NUMBER = {2},
     PAGES = {307--318},
      ISSN = {0021-8693},
   MRCLASS = {20N99 (08A40)},
  MRNUMBER = {682881},
MRREVIEWER = {B.\ H.\ Neumann},
       DOI = {10.1016/0021-8693(82)90305-2},
       URL = {https://doi.org/10.1016/0021-8693(82)90305-2},
}

@article {KS25,
    AUTHOR = {Kaur, Dilpreet and Singh, Pushpendra},
     TITLE = {Classification of primitive quandles of small order},
   JOURNAL = {J. Symbolic Comput.},
  FJOURNAL = {Journal of Symbolic Computation},
    VOLUME = {131},
      YEAR = {2025},
     PAGES = {Paper No. 102452, 9},
      ISSN = {0747-7171,1095-855X},
   MRCLASS = {20N02 (20B15)},
  MRNUMBER = {4897672},
MRREVIEWER = {J.\ D.\ Phillips},
       DOI = {10.1016/j.jsc.2025.102452},
       URL = {https://doi.org/10.1016/j.jsc.2025.102452},
}

@article {Mat82,
    AUTHOR = {Matveev, S. V.},
     TITLE = {Distributive groupoids in knot theory},
   JOURNAL = {Mat. Sb. (N.S.)},
  FJOURNAL = {Matematicheski\u{\i} Sbornik. Novaya Seriya},
    VOLUME = {119(161)},
      YEAR = {1982},
    NUMBER = {1},
     PAGES = {78--88, 160},
      ISSN = {0368-8666},
   MRCLASS = {57M25 (20L15)},
  MRNUMBER = {672410},
MRREVIEWER = {Jonathan\ A.\ Hillman},
}

@article {Pra93,
    AUTHOR = {Praeger, Cheryl E.},
     TITLE = {An {O}'{N}an-{S}cott theorem for finite quasiprimitive
              permutation groups and an application to {$2$}-arc transitive
              graphs},
   JOURNAL = {J. London Math. Soc. (2)},
  FJOURNAL = {Journal of the London Mathematical Society. Second Series},
    VOLUME = {47},
      YEAR = {1993},
    NUMBER = {2},
     PAGES = {227--239},
      ISSN = {0024-6107,1469-7750},
   MRCLASS = {05C25 (20B15)},
  MRNUMBER = {1207945},
MRREVIEWER = {Dragan\ Maru\v{s}i\v{c}},
       DOI = {10.1112/jlms/s2-47.2.227},
       URL = {https://doi.org/10.1112/jlms/s2-47.2.227},
}

@book{Cam99,
  title={Permutation groups},
  author={Cameron, Peter J},
  year={1999},
  publisher={Cambridge University Press}
}

@misc{Ven14,
 title ={rig, Racks in GAP, Version 1},
 author={Vendramin, Leandro},
 month ={May},
 year ={2014},
 note ={GAP package},
 howpublished ={\href{https://gap-packages.github.io/rig/}{\color{black}{\texttt{https://gap-packages.github.io/}\discretionary{}{}{}\texttt{rig/}}}},
 printedkey={Ven14}
}

@misc{Ber22,
title={QuimpGrp, A database of QUasiprimitive IMPrimitive permutation groups},
author={Bernhardt, Dominik},
year ={2022},
note ={GAP package},
howpublished={\href{https://github.com/gap-packages/QuimpGrp}{\color{black}{\texttt{https://github.com/gap-packages/QuimpGrp}}}}
}

@article{GAP,
  title={GAP--Groups, Algorithms, and Programming, Version 4.12.2; 2022},
  author={GAP Group and others},
  journal={URL: https://www.gap-system.org},
  year={2023}
}

@misc{DKa,
 title ={Database of simple quandles},
 author={Kaur, Dilpreet and Singh, Pushpendra},
 howpublished={\href{https://sites.google.com/view/dilpreetkaur/others/database-for-simple-quandles}{\color{black}{\texttt{sites.google.com/\linebreak view/dilpreetkaur/others/database-for-simple-quandles}}}},
 year={2025}

}

@inbook{KingBCC05, place={Cambridge}, series={London Mathematical Society Lecture Note Series}, title={The subgroup structure of finite classical groups in terms of geometric configurations}, booktitle={Surveys in Combinatorics 2005}, publisher={Cambridge University Press}, author={King, Oliver H.}, year={2005}, pages={29–56}, collection={London Mathematical Society Lecture Note Series}}

@article {FGV10,
    AUTHOR = {Freyre, Sebasti\'{a}n and Gra\~{n}a, Mat\'{\i}as and
              Vendramin, Leandro},
     TITLE = {On {N}ichols algebras over {${\rm PGL}(2,q)$} and {${\rm
              PSL}(2,q)$}},
   JOURNAL = {J. Algebra Appl.},
  FJOURNAL = {Journal of Algebra and its Applications},
    VOLUME = {9},
      YEAR = {2010},
    NUMBER = {2},
     PAGES = {195--208},
      ISSN = {0219-4988,1793-6829},
   MRCLASS = {16T05 (20D08 20G40)},
  MRNUMBER = {2646659},
MRREVIEWER = {Istv\'{a}n\ Heckenberger},
       DOI = {10.1142/S0219498810003823},
       URL = {https://doi.org/10.1142/S0219498810003823},
}

@misc{Ada02,
title={Character tables of $\GL(2,q), \SL(2,q), \PGL(2,q)$ and $\PSL(2,q)$},
author={Adams, James},
year ={2002},
howpublished={\linebreak \href{http://www.math.umd.edu/~jda/characters/}{\color{black}{\texttt{http://www.math.umd.edu/~jda/characters/}}}}
}

\Addresses

\end{document}